\documentclass{amsart}
\usepackage{amsmath,amssymb}
\newtheorem{theorem}{Theorem}
\newtheorem{lemma}[theorem]{Lemma}
\newtheorem{proposition}[theorem]{Proposition}
\newtheorem{corollary}[theorem]{Corollary}

\begin{document}
\title{Finitely presented modules over semihereditary rings}
\author{Fran\c cois Couchot}
\address{Laboratoire de Math\'ematiques Nicolas Oresme, CNRS UMR
  6139,
D\'epartement de math\'ematiques et m\'ecanique,
14032 Caen cedex, France}
\email{couchot@math.unicaen.fr} 

\begin{abstract} We prove that each almost local-global semihere\-ditary ring
  $R$ has the stacked bases property and is almost B\'ezout. More
  precisely, if $M$ is a finitely presented module, its torsion part
  $tM$ is a direct sum of cyclic modules where the family of
  annihilators is an ascending chain of invertible ideals. These
  ideals are invariants of $M$. Moreover $M/tM$ is a projective
  module which is isomorphic to a direct sum of finitely generated
  ideals. These ideals allow us to define a finitely generated ideal
  whose isomorphism class is an invariant of $M$. The idempotents
  and the positive integers defined by the rank of $M/tM$ are 
  invariants of $M$ too. It follows that each semihereditary ring
  of Krull-dimension one or of finite character, in particular each hereditary
  ring, has the stacked base property. These results were already
  proved for Pr\"ufer domains by Brewer, Katz, Klinger, Levy and
  Ullery. It is also shown that every semihereditary B\'ezout ring of
  countable character is an elementary divisor ring.
\end{abstract}

\maketitle

\bigskip
It is well-known for a long time that every Dedekind domain satisfies the
\textbf{stacked bases property} (or
the Simultaneous Basis Property). See \cite{Kap52}. The definitions will
be given below. More recently this result was extended to every
Pr\"ufer domain of finite character \cite{Lev87} or of Krull dimension
one \cite{BrKl87} and more generally to each almost local-global
Pr\"ufer domain \cite{BrKl87}.

The aim of this paper is to show that every almost local-global
semihereditary ring, and in particular every hereditary ring, has 
the stacked bases property too. This extension is possible principally
because every semihereditary ring has the following properties:
\begin{itemize}
\item every finitely presented module is
the direct sum of its torsion submodule with a projective module,
\item the annihilator of each finitely generated ideal is generated by
  an idempotent,
\item and each faithful finitely generated ideal contains a
  nonzerodivisor.

\end{itemize}
So, the only difference with the domain case is that we have to manage
nontrivial idempotents. Except for this difference, we do as  
Section 4 of Chapter V of \cite{FuSa01}. 

All rings in this paper are unitary and commutative. 
We say that a ring $R$ is \textbf{of finite (resp. countable)
  character} if every non-zero-divisor is contained in  finitely
(resp. countably) many maximal ideals. We
say that $R$ is {\bf local-global} if each polynomial over $R$ in
finitely many indeterminates which admits unit values locally, admits
unit values. A ring $R$ is said to be {\bf almost local-global} if for
every faithful principal ideal $I$, $R/I$ is local-global. 

Let $R$ be a semihereditary ring. If for each finitely presented
module $M$ and for
each presentation $0\rightarrow H\rightarrow F\rightarrow
M\rightarrow 0$, where $F$ is a finitely generated free module, there exist finitely generated ideals $J_1,\dots,J_n$
and invertible ideals $I_1\supseteq \dots\supseteq I_m$ such that
$F\cong J_1\oplus\dots\oplus J_n$ and $H\cong I_1J_1\oplus\dots\oplus
I_mJ_m$ we say that $R$ has the {\bf stacked bases property}. We will
prove the two following theorems:
\begin{theorem} \label{T:stak} Every almost local-global
  semihereditary ring has the stacked bases property.
\end{theorem}
\begin{theorem} \label{T:decom} Let $R$ be an almost local-global
  semihereditary ring and $M$ a finitely presented module. Then $M$
  decomposes as:
\[M\cong (\bigoplus_{1\leq j\leq m}R/I_j)\oplus\bigl(\bigoplus_{1\leq i\leq
  p}(\bigoplus_{1\leq k\leq n_i}J_{i,k}e_i)\bigr)\]
where $I_1\supseteq\dots\supseteq I_m$ are proper invertible ideals,
$\{e_1,\dots,e_p\}$ is a family of ortho\-gonal idempotents,
$\{n_1,\dots,n_p\}$ a strictly increasing sequence of integers $>0$
and for each $i,\ 1\leq i\leq p$, $J_{i,1},\dots,J_{i,n_i}$ are arbitrary
invertible ideals of $Re_i$. The invariants of $M$ are:
\begin{itemize}
\item[(a)] the ideals $I_1,\dots,I_m$,
\item[(b)] the idempotents $e_1,\dots,e_p$,
\item[(c)] the integers $n_1,\dots,n_p$ and 
\item[(d)] the isomorphism class of \(\oplus_{1\leq i\leq
    p}(\Pi_{1\leq k\leq n_i}J_{i,k})\).
\end{itemize}
($R$ satisfies the Invariant Factor Theorem)
\end{theorem}

\bigskip From these two theorems we deduce the two following corollaries:
\begin{corollary} Let $R$ be a semihereditary ring. Assume that $R$
  satisfies one of the following conditions:
\begin{itemize}
\item[(a)] $R$ is of finite character
\item[(b)] $R$ has Krull-dimension $\leq 1$.
\end{itemize}
Then $R$ has the stacked bases property and every finitely presented
module $M$ has a decomposition as in Theorem~\ref{T:decom}.
\end{corollary}
\textbf{Proof.} Let $I=Rs$ where $s$ is a non-zero-divisor. Then $s$
doesn't belong to any minimal prime ideal. Thus $R/I$ is either
semilocal or zero-Krull-dimensional. So, it is local-global by
\cite[Proposition~p.455]{McWa81} . \qed 

\begin{corollary} Let $R$ be a hereditary ring. Then $R$ has the stacked bases property and every finitely presented
module $M$ has a decomposition as in Theorem~\ref{T:decom}.
\end{corollary}
\textbf{Proof.} $R$ is both one-Krull-dimensional and of finite
character. \qed

\bigskip To prove Theorem~\ref{T:stak} and Theorem~\ref{T:decom} some
preliminary results are needed. It is obvious that each semihereditary
ring is a {\bf pp-ring} (or a {\bf Baer ring}), i.e. a ring for which
every principal ideal is projective.
\begin{lemma} \label{L:quen}
Let $R$ be a pp-ring and $Q$ its ring of quotients. Then the following
assertions hold: 
\begin{enumerate}
\item For each finitely generated ideal $A$, $(0:A)$ is generated by
  an idempotent.
\item each faithful finitely generated ideal
  contains a non-zero-divisor.
\item $Q$ is absolutely flat.
\end{enumerate}
\end{lemma}
\textbf{Proof.} The first assertion is well known if $A$ is
principal. Let $e_1$ and $e_2$ be two idempotents such that
$(0:A_i)=Re_i$, where $A_i$ is an ideal, for $i=1,2$. Then
$(0:A_1+A_2)=Re_1\cap Re_2=Re_1e_2$.

From the first assertion we deduce that for each $a\in R$ there exists
an idempotent 
$e\in R$ such that $ae=0$ and $Ra+Re$ is a faithful ideal. Then the
two last assertions follow from \cite[Proposition~9]{Que71}. \qed

\begin{proposition}
\label{P:tor}
Let $R$ be a semihereditary ring, $M$ a finitely presented $R$-module
and $tM$ its  torsion submodule. 

Then $M\cong tM\oplus M/tM$
and $M/tM$ is projective.
\end{proposition}
\textbf{Proof.} Let $Q$ the ring of fractions of $R$. By
Lemma~\ref{L:quen}, $Q$ is absolutely flat. It follows that
$M/tM\otimes_RQ\cong M\otimes_RQ$ is a projective $Q$-module and
$M/tM$ is flat. By \cite[Proposition~2.3]{CoPe70} ), $M/tM$ is
projective. \qed

\bigskip For an element $x\in R^n$, the {\bf content} $c(x)$ is
defined as the ideal generated by the coordinates of $x$. If $M$ is a
submodule of $R^n$, then $c(M)$ will mean the ideal generated by
$c(x)$ with $x\in M$. If $X$ is a matrix over $R$ then $c(X)$ is the
ideal generated by its entries.

We say that a ring $R$ has the {\bf UCS-property (unit content
  summand)} if, for each $n>0$ the following holds: every finitely
  generated submodule $M$ of $R^n$ with unit content (i.e, $c(M)=R$)
  contains a rank one projective summand of $R^n$.

The following theorem is a generalization of \cite[Theorem~V.4.7]{FuSa01}
and we adapt the proof of \cite[Theorem~3]{BrKl87} to show it.
\begin{theorem} \label{T:UCS} Let $R$ be an almost local-global
  pp-ring. Then $R$ has the UCS-property.
\end{theorem}
\textbf{Proof.} Let $X$ be an $n\times m$-matrix over $R$ such that
$c(X)=R$. Let $M\ (\subseteq R^n)$ be the column space of $X$. There
exist $x_1,\dots,x_p\in M$ such that $c(x_1)+\dots +c(x_p)=R$. For
each $j,\ 1\leq j\leq p$, let $\epsilon_j$ be the idempotent which verifies
$(0:c(x_j))=R(1-\epsilon_j)$. Then for each $j,\ 1\leq j\leq p$, there exists
$a_j\in c(x_j)$ such that $a_1\epsilon_1+\dots+a_p\epsilon_p=1$. By
induction it is 
possible to get an orthogonal family of idempotents $(e_j)_{1\leq
  j\leq p}$ and a family $(b_j)_{1\leq j\leq p}$ such that
$b_1e_1+\dots+b_pe_p=1$ and $e_j\in R\epsilon_j,\ \forall j,\ 1\leq j\leq
p$. After removing the indexes $j$ for which $b_je_j=0$,
we may assume that $b_je_j\not=0,\ \forall j,\ 1\leq j\leq p$. Then
$b_je_j=e_j$, $c(x_j)e_j$ is faithful over $Re_j$ and
$c(e_jX)=c(X)e_j=Re_j$. If we prove that there exists a rank one
summand $U_j$ of $e_jR^n$ contained in $e_jM$ for each $j,\ 1\leq
j\leq p$, then $U_1\oplus\dots\oplus U_p$ is a rank one summand of
$R^n$ contained in $M$. Consequently we may assume that there is only
one idempotent equal to $1$ and $x\in M$ such that $c(x)$ is a
faithful ideal. Then $R^*=R/c(x)$ is a local-global ring. By
\cite[Lemma~V.4.5]{FuSa01} there is a $y^*$ in the image of $M$ such
that $c(y^*)=R^*$. Therefore, for any preimage $y\in M$ of
$y^*,\ c(x)+c(y)=R$ holds.  

Let $N$ be the submodule of $M$ generated by $x$ and $y$. The $2\times
2$-minors of the matrix whose columns are $x$ and $y$ generate an
ideal $I$ in $R$. Let $e$ be the idempotent which verifies
$(0:I)=R(1-e)$. Since $I(1-e)=0$, then $(1-e)N$ is of rank one over
$R(1-e)$, and in view of $c((1-e)x)+c((1-e)y)=R(1-e)$, it is a summand
of $(1-e)R^n$ contained in $(1-e)M$. Since $I$ is faithful over $Re$
then $R^*=Re/I(\cong R/(I+R(1-e))$ is local-global. By
\cite[Lemma~V.4.5]{FuSa01}, $c(x^*+a^*y^*)=R^*$ for some $a^*\in
R^*$. For a preimage $a\in R$ of $a^*$, we have
$c(e(x+ay))+I=Re$. Suppose that $c(e(x+ay))\subseteq P$ for some
maximal ideal $P$ of $Re$. Then $e(x+ay)$ vanishes in $(Re/P)^n$, thus
all the $2\times 2$-minors generators of $I$ are contained in $P$,
i.e., $I\subseteq P$. This is an obvious contradiction to
$c(e(x+ay))+I=Re$. Hence $c(e(x+ay))=Re$, and it follows that the
submodule $T$ generated by $e(x+ay)$ is a rank one summand
of $eR^n$ contained in $eM$. We get that $(1-e)N+T$ is a rank one
 summand of $R^n$ contained in $M$. \qed

\bigskip The following theorem is a generalization of
\cite[Theorem~6]{BKU87} and we do a similar proof. For Pr\"ufer domains
we can also see \cite[Theorem~V.4.8]{FuSa01}.

\begin{theorem} \label{T:sta=ucs} A semihereditary ring $R$ has the
  UCS-property if and only if it has the stacked bases property.
\end{theorem}
\textbf{Proof.}  We do as for Pr\"ufer domains (see
\cite[Theorem~V.4.8]{FuSa01}) to show that ``stacked bases property'' implies
``UCS-property''.

Assume that $R$ has the UCS-property and let $H$ be a finitely
generated submodule of $R^n\ (n\geq 1)$. Let $\epsilon_1$ be the idempotent
such that $(0:c(H))=R(1-\epsilon_1)$. Then $H=\epsilon_1H$ and
$J_1=c(H)\oplus R(1-\epsilon_1)$ is invertible. Thus $J_1^{-1}H$ is a
finitely generated submodule of $\epsilon_1R^n$ and
$c(J_1^{-1}H)=R\epsilon_1$. So, by the UCS-property, $J_1^{-1}H$
contains a rank one summand $U_1$ of $\epsilon_1R^n$, i.e.,
$R^n=U_1\oplus N_1$ for some $N_1$. Hence $J_1^{-1}H=U_1\oplus H_1$
with $H_1=N_1\cap J_1^{-1}H$ and $H=J_1U_1\oplus J_1H_1$. If the
second summand is not $0$, let $\epsilon_2$ be the idempotent such
that $(0:c(H_1))=R(1-\epsilon_2)$ and $J_2=c(H_1)\oplus
R(1-\epsilon_2)$. Then $\epsilon_2\in R\epsilon_1$ and $J_2^{-1}H_1$
contains a rank one summand $U_2$ of $\epsilon_2R^n$. Let $s$ be a
non-zero-divisor such that $sJ_2^{-1}\subseteq R$. If $x\in
J_2^{-1}H_1$ then $x=u+y,\ u\in U_1,\ y\in N_1$ and $sx\in H_1$. It
follows that $su=0$ whence $u=0$. Hence $J_2^{-1}H_1\subseteq N_1$ and 
$U_2$ is a summand of $N_1$. We obtain $R^n=U_1\oplus U_2\oplus N_2$
for some $N_2$. Hence $H=J_1U_1\oplus J_1J_2U_2\oplus J_1J_2H_2$ where
$H_2=N_2\cap J_2^{-1}H_1$. Repetition yields $J_1J_2\dots J_mH_m=0$
\[\mathrm{for}\ m=\mathrm{max}\{\mathrm{rank}(H_P)\mid P\in
\mathrm{Spec}(R)\}\leq n,\ \mathrm{hence} \]
\[R^n=U_1\oplus\dots\oplus U_m\oplus N_m\quad \mathrm{and}\quad H=J_1U_1\oplus\dots\oplus
J_1\dots J_mU_m\]
for rank one projective summand $U_i$ of $\epsilon_iR^n$, where
$(\epsilon_i)_{1\leq i\leq m}$ is a family of nonzero idempotents such
that $\epsilon_{i+1}\in R\epsilon_i$ for every $i,\ 1\leq i\leq
m-1$. Here $N_m$ is projective, so it is isomorphic to a direct sum of
ideals. We thus have stacked bases for $R^n$ and its submodule $H$. \qed

\bigskip
Now it is obvious that Theorem~\ref{T:stak} is an immediate
consequence of Theorem~\ref{T:UCS}~and~\ref{T:sta=ucs}.

\bigskip
We say that a ring $R$ is {\bf almost B\'ezout} if for each faithful
principal ideal $A$, $R/A$ is B\'ezout (i.e., finitely generated
ideals are principal).
The following Proposition~\ref{P:Kap49} was proved by
Kaplansky \cite{Kap52} for Dedekind domains and we do a similar
proof. One can also see
\cite[Proposition~V.4.10]{FuSa01}.

\begin{proposition} \label{P:Kap49} Let $R$ be an almost B\'ezout
  pp-ring. Suppose
\[U=I_1\oplus\dots\oplus I_n,\quad V=J_1\oplus\dots\oplus J_m\]
are decompositions of finitely generated projective $R$-modules $U,\
V$ into direct sums of invertible ideals. Then a necessary and
sufficient condition for $U\cong V$ is that
\[n=m\quad \mathrm{and} \quad I_1\dots I_n\cong J_1\dots J_n.\]
\end{proposition}
\textbf{Proof.} We prove that the condition is necessary as for integral
domains. 

For sufficiency we can also do the same proof since every invertible
ideal $I$ contains a non-zero-divisor $a$. It follows that $aI^{-1}J$
is faithful for each invertible ideal $J$. Therefore
$aI^{-1}/aI^{-1}J$ is a principal ideal of $R/aI^{-1}J$ and generated
by $b+aI^{-1}J$. Then $bR+aI^{-1}J=aI^{-1}$, thus $ba^{-1}I+J=R$. We
  get that $I\oplus J\cong R\oplus IJ$. \qed

\bigskip The following lemma is needed to prove Theorem~\ref{T:decom}.
\begin{lemma} \label{L:ideal}
Let $R$ be an almost local-global semihereditary ring. Then:
\begin{enumerate} 
\item $R$ is an almost B\'ezout ring.
\item For every pair of invertible ideals $I,\ J$, $I/JI\cong R/J$.
\end{enumerate}
\end{lemma}
\textbf{Proof.} Let $J$ be a faithful principal ideal of $R$. Every
finitely generated ideal of $R/J$ is of the form $I/J$ where $I$ is an
invertible ideal of $R$. It follows that $I/JI$ is a rank one
projective $R/J$-module. By \cite[Theorem p.457]{McWa81} or
\cite[Proposition~V.4.4]{FuSa01} $I/JI$ is
cyclic. Thus $R/J$ is a B\'ezout ring: it is also a consequence of
Theorem~\ref{T:ariloglo}. 

Since $JI$ contains a faithful principal ideal, $I/JI$ is cyclic. From
the obvious equality $I(JI:I)=IJ$, we deduce that $(JI:I)=J$. \qed

\bigskip
Now we can prove Theorem~\ref{T:decom}.

\textbf{Proof of Theorem~\ref{T:decom}.} 

Suppose that $M\cong R^n/H$. By
  Theorem~\ref{T:sta=ucs} and its proof we have
\[R^n=U_1\oplus\dots\oplus U_m\oplus N_m\quad \mathrm{and}\quad
H=J_1U_1\oplus\dots\oplus J_1\dots J_mU_m\]
Moreover there is a family of idempotents $(\epsilon_i)_{1\leq i\leq m}$ such
that $U_i$ is isomorphic to an invertible ideal of $R\epsilon_i$ and
$\epsilon_{i+1}\in R\epsilon_i$. We set
$I_k=J_1\dots J_k\epsilon_k\oplus R(1-\epsilon_k)$ for each $k,\ 1\leq
k\leq m$. By Lemma~\ref{L:ideal} 
\[U_i/I_iU_i\cong R\epsilon_i/I_i\epsilon_i\cong R/I_i,\ \forall i,\ 1\leq i\leq
m.\] 
By using \cite[Proposition~1.10]{ShWi74} it is possible to get that
\[tM\cong R/I_1\oplus\dots\oplus R/I_m\ \mathrm{with}\ 
I_1\supseteq\dots\supseteq I_m.\]
 By \cite[Theorem~9.3]{Kap49},
$I_1,\dots,I_m$ are invariants of $tM$ and $M$. 

Since $N_m\cong M/tM$ is a finitely generated projective module
there exists an orthogonal family of idempotents $(e_i)_{1\leq i\leq
  p}$ such that $e_iN_m$ has a constant rank $n_i$ over
  $Re_i$. So $e_iN_m\cong J_{i,1}\oplus\dots\oplus
  J_{i,n_i}$ where $J_{i,k}$ is an invertible ideal of $Re_i$ for each
  $k,\ 1\leq k\leq n_i$. By Proposition~\ref{P:Kap49} the isomorphism
  class of $J_{i,1}\dots J_{i,n_i}$ is an invariant of $e_iN_m$. We
  conclude that the isomorphism class of \(\oplus_{1\leq i\leq
    p}(\Pi_{1\leq k\leq n_i}J_{i,k})\) is an invariant of $M$. \qed

\bigskip It is also possible to deduce Theorem~\ref{T:decom} from the
following. A ring $R$ is said to be \textbf{arithmetic} if it is locally B\'ezout.
\begin{theorem} \label{T:ariloglo} Let $R$ be an arithmetic local-global ring
 . Then $R$ is an elementary divisor ring. Moreover, for each $a,b\in R$, there
  exist $d,a',b',c\in R$ such that $a=a'd$, $b=b'd$ and 
  $a'+cb'$ is a unit of $R$. 
\end{theorem}
\textbf{Proof.}  Since every finitely generated ideal is locally
principal $R$ is B\'ezout by \cite[Corollary 2.7]{EsGu82}. Let $a,b\in
R$. Then there exist $a',b',d\in R$ such that $a=a'd,b=b'd$ and
$Ra+Rb=Rd$. Consider the following polynomial $a'+b'T$. If $P$ is a
maximal ideal, then we have $aR_P=dR_P$ or $bR_P=dR_P$. So, $a'$ or
$b'$ is a unit of $R_P$, whence $a'+b'T$ admits unit values in $R_P$. Then the last assertion holds. Now let $a,b,c\in R$ such that
$Ra+Rb+Rc=R$. We set $Rb+Rc=Rd$. Let $b',c',s$ and $q$ such that
$b=b'd,\ c=c'd$ and $b'+c'q$ and $a+sd$ are units. Then
$(b'+c'q)(a+sd)=(b'+c'q)a+s(b+qc)$ is a unit. We conclude by
\cite[Theorem 6]{GiHe56}.  \qed

\bigskip
\textbf{Another proof of Theorem~\ref{T:decom}.} Let $M$ be a finitely
presented $R$-module. There exists a non-zero-divisor $s\in R$ which
annihilates $tM$. So, $tM$ is a direct sum of cyclic modules since
$R/Rs$ is an elementary divisor ring by Theorem~\ref{T:ariloglo}. \qed

\bigskip Now we give a generalization of \cite[Theorem~III.6.5]{FuSa01}.
\begin{theorem} Every semihereditary B\'ezout ring of countable
  character is an elementary divisor ring.
\end{theorem}
\textbf{Proof.} By \cite[Theorem~5.1]{Kap49} and \cite[Theorem~3.1]{LLS74}
it is sufficient to prove that every $1\ \mathrm{by}\ 2$, $2\
\mathrm{by}\ 1$ and $2\ \mathrm{by}\ 2$ matrix is equivalent to a
diagonal matrix. By \cite[Theorem~2.4]{LLS74} $R$ is Hermite. So it is
enough to prove that the following matrix $A$ is equivalent to a
diagonal matrix:
\[A=\begin{pmatrix}a & 0\\ b & c\end{pmatrix}.\]
Let $e$ the idempotent such that $(0:ac)=R(1-e)$. Thus
$A=eA+(1-e)A$. Since $(1-e)ac=0$ we show, as in the proof of \cite[the proposition]{Sho74}, that there exist two invertible matrices
  $P_1$ and $Q_1$ and a diagonal matrix $D_1$ with entries in $R(1-e)$ such
  that $P_1(1-e)AQ_1=D_1$. Since $Re/Rac\cong R/R(ac+1-e)$ and 
  $(ac+1-e)$ is a non-zerodivisor, $ac$ is contained in countably many
  maximal ideals of $Re$. We prove, as in the proof of \cite[Theorem~III.6.5]{FuSa01}, that there exist two invertible matrices 
  $P_2$ and $Q_2$ and a diagonal matrix $D_2$ with entries in $Re$ such
  that $P_2eAQ_2=D_2$. We set $P=(1-e)P_1+eP_2,\ Q=(1-e)Q_1+eQ_2\
  \mathrm{and}\ D=(1-e)D_1+eD_2$. Then $P$ and $Q$ are invertible, $D$
  is diagonal and $D=PAQ$. \qed

\end{document}